\documentclass[twocolumn,letter]{autart}    % Enable this line and disable the
\everymath{\displaystyle}
\usepackage{graphicx}
\usepackage{amsmath} % assumes amsmath package installed
\usepackage{amssymb}  % assumes amsmath package installed
\usepackage{xcolor}
\usepackage{placeins}
\newcommand{\beqnum}{\begin{equation}\begin{array}{lcl}}
\newcommand{\eeqnum}{\end{array}\end{equation}}
\newcommand{\beq}{\begin{eqnarray*}}
\newcommand{\eeq}{\end{eqnarray*}}
\newcommand{\bef}{\begin{figure}[tbh!]}
\newcommand{\enf}{\end{figure}}

\newtheorem{montheo}{\bf Theorem}

\newtheorem{Proposition}{\bf Proposition}

\newcommand*{\QEDC}{\vspace{-0.4cm}\begin{flushright} $\blacksquare$ \end{flushright}\vspace{-0.4cm}}

\makeatletter
\begin{document}

\begin{frontmatter}
\title{A Lyapunov approach to Robust and Adaptive finite time stabilization of integrator chains with bounded uncertainty
\thanksref{footnoteinfo}}

\thanks[footnoteinfo]{This paper was not presented at any IFAC meeting. Corresponding author S. Laghrouche. Tel: +33 (0)3 84 58 34 19}
%Fax +XXXIX-VI-mmmxxv.
%
\author[set]{Mohamed Harmouche}\ead{\\mohamed.harmouche@utbm.fr},    % Add the
\author[set]{Salah Laghrouche}\ead{salah.laghrouche@utbm.fr},               % e-mail address
\author[lss]{Yacine Chitour}\ead{yacine.chitour@lss.supelec.fr}  % (ead) as shown
\address[set]{SET Laboratory, UTBM, Belfort, France}
\address[lss]{L2S, Universite Paris XI, CNRS 91192 Gif-sur-Yvette, France.}
\begin{abstract}                          % Abstract of not more than 200 words.
In this paper, we present Lyapunov-based robust and adaptive controllers for the finite time stabilization of a perturbed chain of integrators with bounded uncertainties. The proposed controllers can be designed for integrator chains of any arbitrary length. The uncertainty bounds are known in the robust control problem whereas they are unknown in the adaptive control problem. Both controllers are developed from a class of finite time stabilization controllers for pure integrator chains. Lyapunov-based design permits to calculate upper bound on convergence time.
\end{abstract}
\end{frontmatter}
%%%%%%%%%%%%%%%%%%%%%%%%%%%%%%%%%%%%%%%%%%%%%%%%%%%%%%%%%%%%%%%%%%%%%%%%%%%%%%%%
\section{Introduction}
The problem of finite time stabilization of perturbed integrator chains arises in many practical cases of robust nonlinear control. Usually this uncertainty is bounded by the physical limits of the system, however the bounds may be known or unknown.

For the case of known uncertainty bounds, Levant used homogeneity approach to demonstrate finite time stabilization of integrator systems \cite{Levant2003,Levant2005}. Laghrouche et al. \cite{Laghrouche2} presented a two part integral sliding mode based control to deal with the finite time stabilization and uncertainty rejection separately. Dinuzzo et al. proposed another method in \cite{Dinuzzo09}, where finite time stabilization is treated as Robust Fuller's problem using Higher Order Sliding Mode. Defoort et al. \cite{Defoort2009} developed a robust MIMO controller, using a constructive algorithm with geometric homogeneity based finite time stabilization.

The problem is more challenging if the uncertainty bounds are unknown. For this, the control design should (a) not require the uncertainty bounds and (b) avoid gain overestimation \cite{Plestan2010}. Huang et al. \cite{Huang2008} used dynamic gain adaptation for first order systems. Their method does not solve the gain overestimation problem because the gains cannot decrease. Plestan et al. \cite{Plestan2010,Plestan2012} proposed a sliding mode approach, in which the gains are decreased slowly, after sliding mode is reached. Shtessel et al. \cite{Shtessel} also presented a Second Order adaptive gain super-twisting SMC for non-overestimation of the control gains. Glumineau et al. \cite{Glumineau_Adaptive} used impulsive sliding mode for adaptive control of a double integrator system.

In this paper, we present Lyapunov-based controllers for the finite time stabilization of arbitrary order perturbed integrator chains with bounded uncertainties. There are two main contributions in this paper. First, a robust controller is developed that stabilizes an integrator chain of arbitrary length, if the bounds of the uncertainty are known. The advantage is that it can be developed from a class of finite time controllers for pure integrator chains. Then the controller is extended to an adaptive controller for the case where the bounds on the uncertainty are unknown. This controller aims to converge the states to a neighborhood of the origin. However, the states may leave the neighborhood within a region around it, which depends upon the unknown bounds. Therefore, we do not solve the problem of reaching in finite time an arbitrary neighborhood of the origin.

The paper is organized as follows: problem formulation is discussed in Section 2, robust and adaptive controllers are presented in Sections 3 and 4 respectively and conclusion is given in Section 5.
%%%%%%%%%%%%%%%%%%%%%%%%%%%%%%%%%%%%%%%%%%%%%%%%%%%%%%%%%%%%%%%%%%%%%%%%%%%%%%%%
\section{Problem Formulation}
\noindent Let us consider an uncertain nonlinear system:
\vspace{-0.5cm}	\beqnum\label{u.l.s.0}
 \dot z_i &=& z_{i+1}, \ i=1,...,r-1,\\
 \dot z_{r} &=& \varphi(t) + \gamma(t)u .
 \eeqnum
\noindent where $z \in {\mathbb{R}^r}$ is the state vector and $u \in \mathbb{R}$ is the control input. The functions $\varphi$ and $\gamma$ are arbitrary measurable functions that represent bounded uncertainty:
\vspace{-0.3cm}	$$(H1)\quad \varphi(t) \in I_\varphi:=\left[-\bar\varphi , \bar\varphi \right], \quad \gamma(t) \in I_\gamma:= \left[\gamma_m , \gamma_M \right],$$
where $\bar\varphi, \gamma_m , \gamma_M$ are positive constants. In consequence, we are in fact dealing with the differential inclusion 
\begin{equation}\label{u.l.s.}
\dot z_r\in  I_\varphi+u I_\gamma.
\end{equation}
The control objective is to stabilize System \eqref{u.l.s.} to the origin in finite time. Since these controllers are discontinuous feedback laws $u=U(z)$, solutions of  \eqref{u.l.s.} will fall under differential inclusions and need to be understood in Filippov sense, i.e. the right hand vector set is enlarged at the discontinuity points of \eqref{u.l.s.} to the convex hull of the set of velocity vectors obtained by approaching $z$ from all the directions in $\mathbb{R}^r$, while avoiding zero-measure sets \cite{Filippov}.
%--------------------------------------------------------------
\section{Design of robust controller}\label{robust}
\noindent In this section, we develop a controller for stabilizing System \eqref{u.l.s.}, assuming that the bounds on $\varphi$ and $\gamma$ are known. This controller is derived from a class of Lyapunov-based controllers that guarantee finite time stabilization of pure integrator chains, and satisfy some additional geometric conditions. This chain is given as
\vspace{-0.3cm}	\beqnum\label{integrators}
\left\{
\begin{array}{lcl}
 \dot z_i &=& z_{i+1}, \ i=1,...,r-1,\\
 \dot z_{r} &=& u .
 \end{array}
 \right.
\eeqnum
\noindent Let us recall the theorem:
\vspace{-0.3cm}\begin{montheo} \label{Bernstein_theorem}
\cite{Bernstein2000} Consider System \eqref{integrators}. Suppose there exist a continuous state-feedback control law $u=u_0 (z)$, a positive definite $C^1$ function $V_1$ defined on a neighborhood ${{\hat U}} \subset {\mathbb{R}^r}$ of the origin and real numbers $c>0$ and $0< \alpha < 1,$ such that the condition $ \dot{V}_1 + c{V_1}^\alpha(z(t))\leqslant 0, \hbox{ if }z(t) \in {\hat U}$ is true for every trajectory $z$ of System \eqref{integrators}. Then all trajectories of System \eqref{integrators} with the feedback $u_0(z)$ which stay in ${{\hat U}}$ converge to zero in finite time.
If ${{\hat U}} = {\mathbb{R}^r}$ and $V_1$ is radially unbounded, then System \eqref{integrators} with the feedback $u_0(z)$  is globally finite time stable with respect to the origin.
\end{montheo}
\noindent Based on this theorem, we develop a robust controller for System \eqref{u.l.s.}.
\vspace{-0.3cm}\begin{montheo}\label{Harmouche_condition} Consider System \eqref{u.l.s.} subject to Hypothesis \textbf{H1}. Then the following control law stabilizes System \eqref{u.l.s.} to the origin in finite time:
\vspace{-0.5cm}	\beqnum \label{transf}
	u=(u_0 + \bar\varphi sign(u_0))/\gamma_m,
\eeqnum
\noindent where $u_0(z)$ is any state-feedback control law that satisfies the hypotheses of Theorem \ref{Bernstein_theorem} and obeys the following additional conditions: for every $z\in {\hat U}$,
\vspace{-0.3cm}	\beqnum\label{condition_theorem2}
\frac{\partial V_1}{\partial z_r}(z)u_0(z) \leq 0 , \text{and} \,\,\,
{u_0}(z)= 0  \Rightarrow  \frac{{\partial {V_1}}}{{\partial {z_r}}}(z) = 0.
\eeqnum
If ${{\hat U}} = {\mathbb{R}^r}$ and $V_1$ is radially unbounded, then System \eqref{u.l.s.} with the feedback $u(z)$  is globally finite time stable with respect to the origin.
\end{montheo}
\vspace{-0.6cm}\begin{pf*}{Proof of Theorem \ref{Harmouche_condition}.}
\noindent Under the control law $u$ defined in \eqref{transf}, System \eqref{u.l.s.} can be rewritten as:
\vspace{-0.3cm}	\beqnum\label{transf2}
\left\{
\begin{array}{lcl}
 \dot z_i &=& z_{i+1}, \ i=1,...,r-1,\\
 \dot z_{r} &=& \frac{\gamma u_0(z)}{\gamma_m} + \frac{\gamma\bar\varphi}{\gamma_m}sign(u_0(z)) + \varphi .
\end{array}
\right.
\eeqnum
\noindent The Conditions in \eqref{condition_theorem2} mean that $({\partial V_1}/{\partial {z_r}})sign(u_0)$ defines a continuous and non positive function of the time along trajectories of System \eqref{transf2}. The time derivative of the Lyapunov function $V_1$ verifying the hypotheses of Theorem \ref{Bernstein_theorem} along a non trivial trajectory of System \eqref{transf2} inside ${{\hat U}}$ is given as
\vspace{-0.3cm}	\begin{equation*}\begin{array}{lcl}
  \dot V_1 &=& \sum_{i=1}^{r-1} \frac{\partial V_1}{\partial {z_i}}{z_{i+1}} +\frac{{\partial V_1}}{{\partial {z_r}}}\left( {\varphi  + \gamma u} \right),\hfill \\
   &\leqslant& \sum_{i=1}^{r-1} \frac{\partial V_1}{\partial {z_i}}{z_{i+1}}+ \frac{{\partial V_1}}{{\partial {z_r}}}{u_0} + \frac{{\partial V_1}}{{\partial {z_r}}}sign({u_0})\left( {{\bar\varphi} - \left| \varphi  \right|} \right) \hfill \\
  &\leqslant&  \sum_{i=1}^{r-1} \frac{\partial V_1}{\partial {z_i}}{z_{i+1}} + \frac{{\partial V_1}}{{\partial {z_r}}}{u_0} \leq  - c{{V_1}^\alpha } .\hfill
\end{array}\end{equation*}
This implies that any a non-trivial trajectory $z$ reaches zero and stays there in finite time.
\QEDC \end{pf*}
It can be verified that the controllers proposed by Hong \cite{Hong} and Huang \cite{Huang2005} satisfy the hypotheses of Theorem \ref{Bernstein_theorem} and Condition \eqref{condition_theorem2}. Considering Hong's controller, we denote $\lfloor a \rceil^{\theta}:= \left| a \right|^{\theta} sign(a)$ $\forall a\in \mathbb{R}$ and $\theta > 0$. Let $k<0$ and $l_1,\cdots,l_r$ positive real numbers. For $z=(z_1,\cdots,z_r)$, we define for $ i=0,...,r-1$:
\vspace{-0.3cm}	\beqnum\label{Hongfunction}
p_i = 1+(i-1)k,\\ v_0=0,\,
v_{i+1} = -l_{i+1} \lfloor\lfloor z_{i+1} \rceil^{\beta_i } - \lfloor v_i \rceil^{\beta_i } \rceil^{(\alpha_{i+1}/(\beta_i)},
\eeqnum
where $\alpha_i={p_{i+1}}/{p_i}$, for $i=1,...,r$, and, for $k<0$ sufficiently small,
\vspace{-0.3cm}	$$	\beta_0 =p_2,\, (\beta_i + 1)p_{i+1} = \beta_0 + 1 > 0,\,\quad i=1,...,r-1.$$
Consider the positive definite radially unbounded function $V_1:\mathbb{R}^r\rightarrow \mathbb{R}^+$ given by
\vspace{-0.3cm}	\beqnum\label{Lyap-Hong}
{V_1} = \sum\limits_{j = 1}^r \int\limits_{{v_{j - 1}}}^{{z_j}} {{{\left\lfloor s \right\rceil }^{{\beta _{j - 1}}}} - {{\left\lfloor {{v_{j - 1}}} \right\rceil }^{{\beta _{j - 1}}}}ds}.
\eeqnum
It has been proved in \cite{Hong} that for a sufficiently small $k$, there exist $l_i > 0$, $i=1,...,r$, such that the control law $u_0 = v_r$ defined above stabilizes System \eqref{integrators} in finite time and there exists $c>0$ and $0 < \alpha <1$ such that $u_0$ and $V_1$ fulfill the conditions of Theorem \ref{Bernstein_theorem}. Moreover,
\vspace{-0.3cm}	\beqnum
{{\partial V_1}}/{{\partial {z_r}}} &=& {\left\lfloor {{z_r}} \right\rceil ^{{\beta _{r - 1}}}} - {\left\lfloor {{v_{r - 1}}} \right\rceil ^{{\beta _{r - 1}}}},\\
{u_0}(z) &=& {v_r} =  - {l_r}{\left\lfloor {{{\left\lfloor {{z_r}} \right\rceil }^{{\beta _{r - 1}}}} - {{\left\lfloor {{v_{r - 1}}} \right\rceil }^{{\beta _{r - 1}}}}} \right\rceil ^{\frac{\alpha_r}{\beta _{r - 1}}}}.
\eeqnum
It can be verified that
$({\partial V_1}/{\partial z_r}){u_0}(z) \le 0 $
and $u_0(z)=0 \Rightarrow {{\partial V_1}}/{{\partial {z_r}}}=0$.
The feedback law of \cite{Hong} can be simplified by choosing all $\beta_i = 1$ in \eqref{Hongfunction}.
\vspace{-0.3cm}\begin{Proposition}\label{HARMTH}
For System \eqref{integrators}, there exist a sufficiently small  $k<0$ and real numbers $l_i > 0$,
such that the control law $u_0 = v_r$ defined below stabilizes System \eqref{integrators} in finite time. For $i=0,...,r-1$,
\vspace{-0.5cm}	\beqnum\label{Hongsimpl}
v_0=0, \ \  v_{i + 1}=- {l_{i + 1}}{\left\lfloor {{z_{i + 1}} - {v_i}} \right\rceil ^{\frac{1+(i+2)k}{1+(i+1)k}}}.
\eeqnum
\end{Proposition}
\vspace{-1cm}\begin{pf*}{Proof of Proposition \ref{HARMTH}.}
The proof presented in \cite{Hong} is adapted to the parameter choice of \eqref{Hongsimpl}. Let $\lambda =[1+(r+2)k]/[1+(r+1)k]$ and $f_\lambda$ be the closed-loop vector field obtained by using the feedback \eqref{Hongsimpl} in \eqref{integrators}. For each $\lambda > 0$, the vector field $f_\lambda$ is continuous and homogeneous of degree $k < 0$ with respect to the family of dilations $(p_1,...,p_{r})$, where $p_i=1+(i-1)k$, $i=1,...,{r}$. Let $l_i$, $i=1,...,r$ be positive constants such that the polynomial $y^{r} + l_{r}(y^{r-1} + l_{r-1}(y^{r-2} + ... +l_2(y+l_1)))...))$ is Hurwitz. If $k = 0$ the vector field is linear and therefore $\lambda = 1$. Therefore, there exists a positive-definite, radially unbounded, Lyapunov function $V : \mathbb{R}^{r} \rightarrow \mathbb{R}$ such that $L_{f_1}V$ is continuous and negative definite.\\\\
Let $\mathcal{A}=V^{-1}([0,1])$ and $\mathcal{S}=bd \mathcal{A}=V^{-1}(\{1\})$, where $bd \mathcal{A}$ is the boundary of the set $\mathcal{A}$ , i.e. $\mathcal{A}= \left\{z \in \mathbb{R}^r|V(z) \in [0,1]\right \}$  and $\mathcal{S}= \left\{z \in \mathbb{R}^r| V(z) = 1 \right \}$. Then $\mathcal{A}$ and $\mathcal{S}$ are compact since $V$ is proper. Also, $0 \notin \mathcal{S}$ as $V$ is positive definite.
Defining $\phi : (0, 1]\times \mathcal{S} \rightarrow  \mathbb{R}$ by $\phi(\lambda,z)=L_{f_\lambda}V(z)$. Then $V$ is continuous and satisfies $\phi(\lambda,z) < 0$ for all $z \in \mathcal{S} $, i.e. $\varphi(\{1\} \times \mathcal{S}) \subset ( - \infty, 0).$ Since $ \mathcal{S} $ is compact, by continuity there exists $ \epsilon > 0$ such that $\phi((1 - \epsilon, 1] \times {\mathcal{S}}) \subset (- \infty, 0)$. It follows that for $\lambda \in (1 - \epsilon, 1]$, $L_{f_\lambda}V$ takes negative values on $\mathcal{S}$. Thus, $\mathcal{A}$ is strictly positively invariant under $f_\lambda$ for every $\lambda \in (1 - \epsilon, 1]$. Therefore the origin is global asymptotic stable under $f_\lambda$, for $ \lambda  \in (1-\epsilon,1] $.
Finally, for $ \lambda  \in (1-\epsilon,1) $ i.e  $|k|$ small enough, by homogeneity, the origin is globally finite time stable.
\QEDC \end{pf*}
%--------------------------------------------------------------
\section{Adaptive Controller}\label{adaptatif}
\noindent Let us now consider that uncertainty bounds $\gamma_m$, $\gamma_M$ and $\bar \varphi$ of System \eqref{u.l.s.} are unknown. For any $a \in \mathbb{R}$, let $\sigma(a)$ be the standard saturation function defined by $\sigma (a) = {a}/{\max (1, \left| a \right| )}$.
For $\varepsilon > 0$, $a \in \mathbb{R}$, we define $\nu _\varepsilon (a)= 0.5 + 0.5\sigma \left( \left( |a|  - 0.75 \varepsilon \right) /(0.25\varepsilon)  \right)$.
The following controller is proposed:
\vspace{-0.5cm}	\beqnum \label{eq_v}
u = \hat \gamma u_0(z) + \hat \varphi sign( u_0(z) ),
\eeqnum
where  $u_0$ is a homogeneous controller that satisfies the hypotheses of Theorem \ref{Bernstein_theorem} and fulfills Condition \eqref{condition_theorem2}. The adaptive function $\hat \gamma = \kappa + \delta |u_0(z)|$ and  $\hat \varphi(t)$ is defined by the ODE $\dot {\hat \varphi}(t) = {k}{\nu _\varepsilon }({V_1(z)}) - \left( {1 - {\nu _\varepsilon }({V_1(z)})} \right){\left\lfloor {\hat \varphi} \right\rceil ^\eta }$, with the initial condition $\hat \varphi(0) = 0$. The new terms are defined as $\kappa,\ \delta > 0$, $\eta \in (0,1), \ k > 0 $ and $V_1$ is a homogeneous Lyapunov function which also satisfies Theorem \ref{Bernstein_theorem} and Condition \eqref{condition_theorem2}. Then the following theorem provides the main result for the adaptive case.
%%%%%%%%%%%%%%%%%%%%%%%%%%%%%%%%%%%%%%%%%%%%%%%%%%%%%%%%%%%%%%%%%%%%%%%%%%%%%%%%%%%%%%%%%%%%%%%%%%%%%%%%%%%%%%%%%%%%%%%%%%%%%
\vspace{-0.2cm}\begin{montheo}\label{theorem_adaptatif}
Consider System \eqref{u.l.s.} under the feedback control law \eqref{eq_v}. Then, $\forall \varepsilon,\,\,\,\exists \Delta, c'>0$ and $0<\alpha'<1$ such that the following conditions are satisfied for any initial condition $z_0 \in \hat U$
\vspace{-0.3cm}	\begin{description}
\item[$(i)$] $\liminf_{t\rightarrow \infty}V_1(z(t))\le \varepsilon$, $\limsup_{t\rightarrow \infty}V_1(z(t))\le \Delta$;
\item[$(ii)$] $\limsup_{t\rightarrow \infty}{\left| \hat \varphi \right|}\leq 2\bar \Phi + k \left (\Delta^{1 - \alpha}/(c(1 - \alpha))\right),$
\end{description}
where
\vspace{-0.5cm}	$$\bar \Phi := \left( \bar \varphi + {\left(\kappa \gamma_m -1\right)^2}/({4\gamma_m \delta}) \right)/{\gamma_m},$$
\vspace{-0.5cm}	$$\Delta := \left( \varepsilon^{ 1 - \alpha'} + {c'(1-\alpha') \gamma_m {\bar\Phi}^2}/(2k)\right)^{\frac1{1 - \alpha'}}.$$
\end{montheo}
\vspace{-0.6cm}\begin{pf*}{Proof of Theorem 3:} We first demonstrate that when the system states are in the domain $V_1 > \varepsilon$, the controller brings them to the domain $V_1 \le \varepsilon$ in finite time. Then, it is proved that once $z$ reaches the domain {$V_1 \leq \varepsilon$}, it stays in the domain $V_1 \le \Delta$ for all consecutive time instances and $\hat \varphi$ is upper-bounded after a sufficiently large time. It can be noted that $\hat \varphi$ is a non-negative function.

We argue by contradiction in order to prove that $\liminf_{t\rightarrow \infty}V_1(z(t))\le \varepsilon$. Supposing there exists $\bar t$ such that $V_1(t)> \varepsilon$ for every $t\geq \bar t$, then according to the dynamics of ${\hat \varphi}$ , we get $\dot {\hat \varphi} = k$ for $t \ge \bar t$. This implies that for $t \ge \bar t$, $\hat \varphi$ is increasing and $\hat \varphi > \bar\Phi$. Since we have
\vspace{-0.3cm}	\beqnum
	\dot V_1 &=&  \frac{\partial V_1}{\partial z_1 } z_2 + ... + \frac{\partial V_1}{\partial z_r } \left( \gamma  \left[ \hat \gamma u_0 + \hat \varphi sign (u_0) \right] + \varphi \right)\\
		&=& \frac{\partial V_1}{\partial z_1 } z_2 + ... + \frac{\partial V_1}{\partial z_r } u_0 \\
&&+ \frac{\partial V_1}{\partial z_r }\left( -u_0 \!+\! \kappa \gamma u_0 \!+\! \gamma \delta \left\lfloor u_0 \right\rceil ^2  \!+\! \gamma \hat \varphi sign (u_0)  \!+\!  \varphi \right),\\
	&\le& -c V_1 ^\alpha \!-\! \!\left|\frac{\partial V_1}{\partial z_r } \right| \!\! \left( \! (\kappa \gamma_m \!-\! 1) |u_0| \!+\! \gamma_m \delta |u_0|^2  \!+\! \gamma_m \hat \varphi \!-\! \bar \varphi \right),\\
&=&-c V_1 ^\alpha - \left|\frac{\partial V_1}{\partial z_r } \right|
			\left[ \gamma_m \delta \left( |u_0| +\frac{\kappa \gamma_m -1}{2\gamma_m \delta}  \right)^2 \right]\\
&&- \left|\frac{\partial V_1}{\partial z_r } \right|\left[ - \left( \bar \varphi + \frac{\left( \kappa \gamma _m -1\right)^2}{4\gamma_m \delta} \right) + \gamma_m \hat \varphi \right],\\
&\le& -c V_1 ^\alpha - \gamma_m\left|\frac{\partial V_1}{\partial z_r } \right| \left( \hat \varphi - \bar \Phi\right) \le -c V_1 ^\alpha.	
\eeqnum
%which can be rewritten as
%\beqnum
%\dot V_1&\le& -c V_1 ^\alpha - \gamma_m\left|\frac{\partial V_1}{\partial z_r } \right| \left( \hat \varphi - \bar \Phi\right) \le -c V_1 ^\alpha,
%\eeqnum
Then $V_1(z)$ converges to zero in finite time, which contradicts the hypothesis. The functions $u_0$ and $V_1$ are homogeneous, which according to \cite{Bernstein2005}, means that
\vspace{-0.3cm}\beqnum\label{c_a_bernstein}
	\exists \,\,\, c',\,\alpha'>0:  \left|{\partial V_1}/{\partial {z_r}} \right| \le c' {V_1}^{\alpha'},
\eeqnum
where
$ c' = \max_{\{ z:V_1(z) = 1 \}}{ \left| {\partial V_1}/{\partial {z_r}} \right| },\,\,\,  \alpha' = {\kappa_2}/{\kappa_1}$.
The terms $\kappa_2$ and $\kappa_1$ are the respective degrees of homogeneity of ${\partial V_1}/{\partial {z_r}}$ and $V_1$. We suppose now that {$V_1 < \varepsilon$}. Considering \eqref{c_a_bernstein}, let us estimate the overshoot in the worst case condition with respect to uncertainty. For $V_1(z(0)) = \varepsilon$ and $\hat \varphi(0) = 0$, we get
\vspace{-0.3cm}	\beqnum
	\dot V_1 &\leq& -c V_1^\alpha - \gamma_m c' V_1^{\alpha'} \left( \hat \varphi - \bar \Phi \right),\,\, \dot {\hat \varphi} = k.
\eeqnum
 The overshoot $\Delta$ of $V_1$ holds for $\dot V_1 = 0$ at $t=T_M$. We get $\hat \varphi(T_M) = \bar \Phi -\left(\left(c\Delta^{\alpha-\alpha'}\right)/\left(c'\gamma_m\right)\right)  \le \bar \Phi$, and then $T_M \le \bar \Phi/k$. An upper bound of $\Delta$ can be estimated as $\Delta\! =\! \left( \varepsilon^{ 1 - \alpha'} + \left(c'(1-\alpha') \gamma_m {\bar\Phi}^2/(2k)\right)\right)^{\frac1{1 - \alpha'}}$.\\
We now estimate an upper bound of $\limsup_{t\rightarrow \infty}{ \hat \varphi}$. Consider the case $V_1(z(0)) = \varepsilon$ with $\dot V_1(z(0)) \ge 0$, in this case we have $\hat \varphi (0) < \bar \Phi$. For $t = T_M$, i.e., $\dot V_1 = 0$, we get $\hat \varphi (T_M) \le \bar \Phi + \hat \varphi (0) \le 2\bar \Phi$. $\hat \varphi$ will increase until time $T_f$ where $\dot{\hat \varphi}(T_f) = 0$ and $V_1(z(T_f)) \ge 0$.
The worst case is calculated with respect to the boundary of $\hat \varphi$, using $\dot V_1  \leq -c V_1^{\alpha}$ and $\dot{\hat \varphi} = k$. Here $T_f$ corresponds to $V_1(z(T_f)) = 0$, i.e $T_f - T_m = \left(\Delta^{(1-\alpha)}\right)/\left(c(1-\alpha)\right)$, which implies that
\vspace{-0.5cm}	\beqnum\hat \varphi(T_f) &\leq& \hat \varphi(T_M) + k(T_f - T_M) \\&=& 2 \bar \Phi + \left(k \Delta^{(1-\alpha)}\right)/\left(c(1-\alpha)\right).\eeqnum
\vspace{-0.5cm}\QEDC \end{pf*}
\vspace{-0.3cm}\textbf{Discussion:} The adaptive functions $\hat \gamma$ and $\hat \varphi$ are chosen non-negative and $\hat \gamma$ is strictly positive with a term proportional to $| u_0 |$. It satisfies the condition  ${\partial \hat \gamma}/{\partial | u_0 | } >0$, which is sufficient to ensure that the states will not diverge, irrespectively of  $\hat \varphi$. The second adaptive function $\hat \varphi$ ensures the convergence of the state to a neighborhood of zero. Its dynamics can be defined explicitly by:
\vspace{-0.5cm}	\beqnum \label{dphi_hat}
\dot { \hat \varphi} = \left  \{ \!{\begin{array}{{ccc}}
{{k} }&,&{{V_1} \ge \varepsilon }\\
\left(\! V_1\!-\!\frac{\varepsilon}2\! \right)\frac{2k}{\varepsilon}\! - \!\left(\varepsilon-V_1\right)\frac{2}{\varepsilon}{{\left\lfloor {\hat \varphi} \right\rceil }^\eta }\! \!&,&\!\frac{\varepsilon }{2}\! \le\! {V_1}\! \le\! \varepsilon, \\
{ - {{\left\lfloor {\hat \varphi} \right\rceil }^\eta }  } &,&{{V_1}} \le \frac{\varepsilon }{2} .
\end{array}} \right.
\eeqnum
\vspace{-0.2cm}\section{CONCLUSIONS}\vspace{-0.2cm}
\noindent In this paper, we have presented a robust and an adaptive controller for the finite time stabilization of perturbed integrated chains with bounded uncertainty. The robust controller is designed using the knowledge of the uncertainty bounds, and it converges exactly to zero. In the case of unkwown uncertainty bounds, the adaptive controller make the state enter in finite time an arbitrary neighborhood of zero an infinite number of time while staying in a bounded neighborhood of zero of depending of the uncertainty bounds. The proof of convergence of both controllers has been demonstrated through Lyapunov analysis, and the calculation of upper bound of the convergence time has also been presented.
%%%%%%%%%%%%%%%%%%%%%%%%%%%%%%%%%%%%%%%%%%%%%%%%%%%%%%%%%%%%%%%%%%%%%%%%%%%%%%%%

%\FloatBarrier
%\clearpage
%\newpage
\footnotesize{\bibliography{References_Stabilisation_Automatica2012}}

\begin{thebibliography}{10}

\bibitem{Levant2003}
A.~Levant.
\newblock Higher-order sliding modes, differentiation and output-feedback
  control.
\newblock {\em International Journal of Control}, 76(9/10):924 -- 941, 2003.

\bibitem{Levant2005}
A.~Levant.
\newblock Homogeneity approach to high-order sliding mode design.
\newblock {\em Automatica}, 41(5):823 -- 830, 2005.

\bibitem{Laghrouche2}
S.~Laghrouche, F.~Plestan, and A.~Glumineau.
\newblock Higher order sliding mode control based on integral sliding mode.
\newblock {\em Automatica}, 43(3):531--537, 2007.

\bibitem{Dinuzzo09}
F.~Dinuzzo and A.~Fererra.
\newblock {Higher Order Sliding Mode Controllers with Optimal Reaching}.
\newblock {\em IEEE Transactions on Automatic Control}, 54(9):2126--2136, 2009.

\bibitem{Defoort2009}
M.~Defoort, T.~Floquet, A.~Kokosy, and W.~Perriquetti.
\newblock A novel higher order sliding mode control scheme.
\newblock {\em Systems and Control Letters}, 58(2):102 -- 108, 2009.

\bibitem{Plestan2010}
F.~Plestan, Y.~Shtessel, V.~Brégeault, and A.~Poznyak.
\newblock New methodologies for adaptive sliding mode control.
\newblock {\em International Journal of Control}, 83(9):1907 -- 1919, 2010.

\bibitem{Huang2008}
Y.J. Huang, T.C. Kuo, and S.H. Chang.
\newblock Adaptive sliding-mode control for nonlinear systems with uncertain
  parameters.
\newblock {\em IEEE Transactions on Automatic Control}, 38(2):534--539, 2008.

\bibitem{Plestan2012}
F.~Plestan, Y.~Shtessel, V.~Br\'egeault, and A.~Poznyak.
\newblock Sliding mode control with gain adaptation - application to an
  electropneumatic actuator (in press).
\newblock {\em Control Engineering Practice}, 2012.
\newblock doi: 10.1016/j.conengprac.2012.04.012.

\bibitem{Shtessel}
Y.~Shtessel, M.~Taleb, and F.~Plestan.
\newblock A novel adaptive gain supertwisting sliding mode controller:
  methodology and application.
\newblock {\em Automatica}, 48(5):759--769, 2011.

\bibitem{Glumineau_Adaptive}
A.~Glumineau, Y.~Shtessel, and F.~Plestan.
\newblock Impulsive sliding mode adaptive control of second order system.
\newblock pages 5389 -- 5394, Milano, Italy, August - September 2011.

\bibitem{Filippov}
A.F. Filippov.
\newblock {\em Differential Equations with Discontinuous Right-Hand Side}.
\newblock Kluwer, Dordrecht, The Netherlands, 1988.

\bibitem{Bernstein2000}
S.P. Bhat and D.S. Bernstein.
\newblock Finite-time stability of continuous autonomous systems.
\newblock {\em SIAM Journal of Control and Optimization}, 38(3):751 -- 766,
  2000.

\bibitem{Hong}
Y.~Hong.
\newblock Finite-time stabilization and stabilizability of a class of
  controllable systems.
\newblock {\em Systems and Control Letters}, 46(4):231--236, 2002.

\bibitem{Huang2005}
X.~Huang, W.~Lin, and B.~Yang.
\newblock Global finite-time stabilization of a class of uncertain nonlinear
  systems.
\newblock {\em Automatica}, 41(5):881--888, 2005.

\bibitem{Bernstein2005}
Sanjay~P Bhat and Dennis~S Bernstein.
\newblock Geometric homogeneity with applications to finite-time stability.
\newblock {\em Math. Control Signals Systems}, 17:101 -- 127, 2005.

\end{thebibliography}
\end{document}